# ROBUST NONPARAMETRIC INFERENCE FOR THE MEDIAN

By Víctor J. Yohai[1] and Ruben H. Zamar[2]

*University of Buenos Aires and CONICET, and University of British Columbia*

We consider the problem of constructing robust nonparametric confidence intervals and tests of hypothesis for the median when the data distribution is unknown and the data may contain a small fraction of contamination. We propose a modification of the sign test (and its associated confidence interval) which attains the nominal significance level (probability coverage) for any distribution in the contamination neighborhood of a continuous distribution. We also define some measures of robustness and efficiency under contamination for confidence intervals and tests. These measures are computed for the proposed procedures.

**1. Introduction.** Often, a fraction of the data is contaminated by outliers and other type of low quality observations. For example, a slight shift in one of several similar instruments used in an experiment may cause a small but consistent bias in a few observations. We are often interested in drawing inference from the uncontaminated part of the data, which distribution we call the "*target distribution*." It is well known that robust point estimates successfully limit the effect of a small fraction of contamination in the data. Unfortunately, naive "robust" confidence intervals constructed around robust point estimates are not that successful. See Fraiman, Yohai and Zamar (2001).

To allow for a fraction $\varepsilon$ of contamination in the data we assume that the actual distribution $G$ belongs to the *contamination neighborhood* of the target distribution $F$,

(1.1) $$\mathcal{F}_\varepsilon(F) = \{G : G = (1-\varepsilon)F + \varepsilon H\},$$

Received February 2003; revised March 2004.
[1]Supported in part by a grant from the Fundación Antorchas, Argentina, Grant X611 from University of Buenos Aires and Grant 03-06277 from the Agencia Nacional de Promoción Científica y Tecnológica, Argentina.
[2]Supported in part by an NSERC operating grant.
*AMS 2000 subject classifications.* Primary 62F35; secondary 62G35.
*Key words and phrases.* Confidence interval, nonparametric, robust, two-sided test.

where $H$ is arbitrary and $0 \leq \varepsilon < 1/2$.

Robust inference (beyond point estimation) means that the inference procedure achieves its intended goal over the entire contamination neighborhood. For instance, robust confidence intervals must achieve the nominal coverage probability of the target parameter for all the distributions in a contamination neighborhood. Similarly, the rejection probability of robust tests when the null hypothesis is true must be smaller than or equal to the nominal significance level under all the distributions in the neighborhood.

Robust tests and confidence intervals have been proposed and studied by several researchers. Huber (1965) introduced censored likelihood ratio tests to robustify the Neyman–Pearson optimal test. Huber (1968) considered robust confidence intervals for a location parameter $\theta$ which cover the true parameter with the nominal probability for all distributions in a neighborhood of the target distribution. The intervals are of the form $(T_n - a, T_n + a)$, where $T_n$ is a location estimate. He found the estimate $T_n$ that minimizes $a$ subject to the conditions $P(T_n < \theta - a) \leq \alpha/2$ and $P(T_n > \theta + a) \leq \alpha/2$—instead of the more natural but less tractable condition $P(T_n < \theta - a) + P(T_n > \theta + a) \leq \alpha$—for finite samples coming from distributions in the contamination neighborhood. The optimal estimate is an M-estimate with Huber type score function. In Huber's approach the scale parameter is assumed known. Fraiman, Yohai and Zamar (2001) solved a related problem: find robust intervals $(T - a, T + a)$ of minimum length and asymptotically correct coverage for all distributions in a contamination neighborhood.

We now briefly discuss two asymptotic approaches to the problem of robust inference for the case of small $\varepsilon$. The first, introduced by Huber-Carol (1970), Rieder (1978) and Bednarski (1982), uses shrinking contamination neighborhoods (contamination fraction of order $n^{-1/2}$) for the null hypothesis and contiguous alternatives of order $n^{-1/2}$. The second, introduced by Rousseeuw and Ronchetti (1981), is based on the influence function for tests which is used to approximate the maximum level and the minimum power of a test in a contamination neighborhood of size $\varepsilon$, when $\varepsilon$ is small. In particular, the approximation of the maximum level can be used to correct the test so that the maximum level is not larger than a given value $\alpha$ for all distributions in a contamination neighborhood. For a full account of this approach see Hampel et al. (1986) and Markatou and Ronchetti (1997). A related approach was given by Lambert (1981) who defines an influence function that measures the effect of the contamination on the $p$-value of a test.

Morgenthaler (1986) considers a class of robust confidence intervals, called strong confidence intervals, which keep the nominal coverage probability conditional on the sample configuration, under two or more specified symmetric distributions. It would seem reasonable to expect that by choosing some extreme symmetric distributions (e.g., the normal and slash distributions), the coverage of the interval should remain correct for other "intermediate" symmetric distributions. Morgenthaler also considers a class of robust confidence intervals, called bioptimal, which are robust in terms of efficiency for two symmetric distributions. The case of asymmetric contamination is not considered in Morgenthaler's approach.

Rieder (1982) addresses the problem of robustifying rank tests preserving their nonparametric nature. He considers one-sided tests for one and two sample problems, showing that the least favorable distribution under a given fraction of contamination does not depend on the target model. Our two-sided modified sign test and the corresponding robust confidence interval can be considered extensions of Rieder's approach.

The rest of the paper is organized as follows. Section 1.1 briefly reviews nonparametric intervals obtained by inverting the sign test. Section 1.2 contains our main result, Theorem 1, which shows that sign-test intervals are not robust and paves the way for the construction of robust nonparametric intervals for the median in Section 2. In this section we also discuss coverage robustness of confidence intervals and the associated concept of level robustness of a test. In Section 3 we address the concept of length robustness of a confidence interval and the associated concept power robustness of a test. In this section we show that the nonparametric robust confidence interval defined in Section 2 has optimal length. In Section 4 we discuss possible extensions and further research. The last section is the Appendix with some proofs. Detailed proofs of our results can be found in Yohai and Zamar (2004).

1.1. *Robust nonparametric inference for the median.* Let

$$x_{(1)} \leq x_{(2)} \leq \cdots \leq x_{(n)}$$

be the order statistics of a sample $X_n = (x_1, \ldots, x_n)$ with common distribution $F$ satisfying the following assumption.

(A1) $F$ is continuous with a unique median $\theta(F) = F^{-1}(1/2)$.

Consider the null hypothesis $H_0: \theta = \theta_0$ and the sign test statistic

$$(1.2) \qquad T_{n,\theta}(X_n) = \sum_{i=1}^{n} I(x_i - \theta > 0).$$

The interval

$$(1.3) \qquad I_\alpha(X_n) = [x_{(k+1)}, x_{(n-k)})$$

4    V. J. YOHAI AND R. H. ZAMARis obtained by inverting the acceptance region $k < T_{n,\theta_0}(X_n) < n - k$. See, for instance, Hettmansperger (1984). The interval (1.3) is a distribution-free $(1 - \alpha(k))100\%$ confidence interval for $\theta$, where

(1.4) $\qquad \alpha(k) = 2P(Z_n \leq k), \qquad Z_n \sim \text{Binomial}(n, 1/2).$

For simplicity, we will only consider levels in the set $\{\alpha(k)\}$, $k = 1, 2, \ldots, [n/2]$. Hettmansperger and Sheather (1986) show how general levels can be obtained by interpolating between the order statistics.

Interval (1.3) yields valid inference for the median of the contaminated distribution, but not for the median of the target distribution. In general, distribution-free methods do not yield valid inference for the target distribution in the presence of asymmetric contamination. Since the median is a very robust location parameter, $\theta(G)$ and $\theta(F)$ are generally close for all $G$ in $\mathcal{F}_\varepsilon(F)$. Still, as shown by Table 1 computed using the result of Theorem 1, the probability that (1.3) covers the target median $\theta(F)$—and the significance level of the associated sign test—may be severely upset.

1.2. *Our main result.* Theorem 1 shows that the nonparametric interval (1.3) is not robust because its probability of covering the median of $F$ can be much smaller than $1 - \alpha(k)$ for distributions $G$ in $\mathcal{F}_\varepsilon(F)$. More importantly, it gives a simple way to modify the definition of this interval (see Section 2.2) so that it remains nonparametric and achieves robustness.

THEOREM 1. *Let $X_n = (x_1, \ldots, x_n)$ be a random sample from $G \in \mathcal{F}_\varepsilon(F)$ with $F$ satisfying* (A1). *Then,*

(a)

(1.5) $\qquad \inf_{G \in \mathcal{F}_\varepsilon(F)} P_G(x_{(k+1)} \leq \theta < x_{(n-k)}) = 1 - \alpha^*(n, k, \varepsilon),$

TABLE 1
*Minimum coverage probability for contaminated samples*

|  | | $\varepsilon$ | | | | | | |
|---|---|---|---|---|---|---|---|---|
| | $1 - \alpha \approx 0.95$ | | | | $1 - \alpha \approx 0.90$ | | | |
| $n$ | 0 | 0.05 | 0.10 | 0.15 | 0 | 0.05 | 0.10 | 0.15 |
| 20 | 0.959 | 0.954 | 0.938 | 0.912 | 0.885 | 0.876 | 0.849 | 0.804 |
| 40 | 0.962 | 0.952 | 0.922 | 0.868 | 0.919 | 0.904 | 0.859 | 0.784 |
| 100 | 0.943 | 0.912 | 0.815 | 0.655 | 0.911 | 0.872 | 0.755 | 0.578 |
| 200 | 0.944 | 0.881 | 0.689 | 0.414 | 0.896 | 0.811 | 0.582 | 0.307 |
| 500 | 0.946 | 0.789 | 0.376 | 0.074 | 0.902 | 0.702 | 0.279 | 0.043 |
| 1000 | 0.946 | 0.636 | 0.108 | 0.002 | 0.906 | 0.537 | 0.068 | 0.001 |
| 2000 | 0.948 | 0.385 | 0.006 | 0 | 0.897 | 0.273 | 0.002 | 0 |

*where*

(1.6) $$\alpha^*(n, k, \varepsilon) = 1 - P(k < Z_n < n - k),$$

*with $Z_n$ distributed as* Binomial$\{n, (1-\varepsilon)/2\}$.

(b) *The infimum in* (1.5) *is achieved for any contaminating distribution which places all its mass to the right or left of $\theta$.*

Using Theorem 1, we calculate the minimum coverage probability for the intervals (1.3) for several values of $n$, $\alpha$ and $\varepsilon$. The results shown in Table 1 are disappointingly low, especially for large $n$. The minimum coverages are not overly pessimistic since they are caused by any contamination fully supported to the right (or left) of the target median.

## 2. Coverage and level robustness.

2.1. *Coverage robustness of a confidence interval.* In connection with the preceding discussion, we now formally state the desired robustness and nonparametric properties for the coverage probability of confidence intervals.

DEFINITION I1 (Coverage robustness). We say that a confidence interval $I_n = [a_n(X_n), b_n(X_n))$ has $\varepsilon$-robust coverage $1 - \alpha$ at $F$ if

(2.1) $$\inf_{G \in \mathcal{F}_\varepsilon(F)} P_G\{a_n(X_n) \leq \theta < b_n(X_n)\} = 1 - \alpha.$$

A related concept of robust confidence interval was introduced by Huber (1968). Although Huber's objective function is not exactly equal to the minimum coverage probability, it is closely related to it. The following definition seems natural to convey the nonparametric nature of an interval.

DEFINITION I2 (Nonparametric coverage robustness). We say that a confidence interval $I_n = [a_n(X_n), b_n(X_n))$ has nonparametric $\varepsilon$-robust coverage $1 - \alpha$ if it has $\varepsilon$-robust level $1 - \alpha$ at $F$ for all $F$ satisfying (A1).

2.2. *An exact nonparametric $\varepsilon$-robust interval for $\theta$.* We wish to construct robust and nonparametric confidence intervals for the median of the target distribution. Theorem 1 derives the exact finite sample least favorable distribution (under contamination neighborhoods) for (1.3) and shows that this distribution does not depend on the target distribution $F$. This theorem also tells us how to modify the interval (1.3) so that it attains nonparametric $\varepsilon$-robust level $1 - \alpha$. Namely, the integer $k$ must satisfy the equation

(2.2) $$\alpha^*(n, k, \varepsilon) = \alpha.$$

Note that the definition (2.2) of $k$ is based on the distribution Binomial$\{n, (1-\varepsilon)/2\}$ instead of the Binomial$(n, 1/2)$. As in the classical case, it is not possible to achieve all the desired exact coverage probabilities $1 - \alpha$. For simplicity, we restrict attention to integers

$$(2.3) \qquad k_n = k_n(n, \alpha) = \arg\min |\alpha^*(n, k, \varepsilon) - \alpha|,$$

which clearly satisfies

$$\lim_{n \to \infty} \alpha^*(n, k_n, \varepsilon) = \alpha.$$

In summary, the modified interval covers the median of the target distribution with a guaranteed confidence level for each $n$ and for all the distributions in a contamination neighborhood of a general target distribution.

2.3. *Level robustness of a test.* Given the well-known duality between confidence interval and tests, it is natural to expect that the nonparametric robust confidence intervals introduced in the previous section will automatically yield nonparametric tests with good robustness properties.

Following Huber (1965), we next define the concept of *$\varepsilon$-robust level-$\alpha$* test.

DEFINITION T1 (Level robustness). Let $F$ be a fixed distribution satisfying (A1) with $\theta = \theta_0$. A nonrandomized test $\varphi_{\theta_0}$ has $\varepsilon$-robust level $\alpha$ (for $H_0$ versus $H_1$) at $F$ if

$$\sup_{G \in \mathcal{F}_\varepsilon(F)} P_G\{\varphi_{\theta_0}(X_n) = 1\} = \alpha.$$

This property ensures the validity of the test over the entire neighborhood $\mathcal{F}_\varepsilon(F)$. That is, the probability of rejecting $H_0$ is less than or equal to $\alpha$ not only at $F$, but also at any $G$ in $\mathcal{F}_\varepsilon(F)$.

DEFINITION T2 (Nonparametric level robustness). We say that a nonrandomized test $\varphi_{\theta_0}$ has nonparametric $\varepsilon$-robust level $\alpha$ (for $H_0$ versus $H_1$) if $\varphi_{\theta_0}$ has $\varepsilon$-robust level $\alpha$ at $F$ for all $F$ satisfying (A1) with $\theta = \theta_0$.

2.4. *An exact nonparametric $\varepsilon$-robust test.* It is immediate that T1 (T2) holds for a family of tests if and only if I1 (I2) holds for the associated sequence of intervals. In particular, the $\varepsilon$-robust sign test $\varphi_{\theta_0}$ of level $\alpha$ can be derived from the nonparametric $\varepsilon$-robust interval $I_\alpha(X_n)$ as follows:

$$\varphi_{\theta_0}(X_n) = \begin{cases} 1, & \text{if } \theta_0 \notin I_\alpha(X_n), \\ 0, & \text{if } \theta_0 \in I_\alpha(X_n), \end{cases}$$

and, therefore,

$$(2.4) \qquad \varphi_{\theta_0}(X_n) = \begin{cases} 1, & \text{if } T_{n,\theta_0}(X_n) \leq k \text{ or } T_{n,\theta_0}(X_n) \geq n - k, \\ 0, & \text{if } k < T_{n,\theta_0}(X_n) < n - k, \end{cases}$$

where $T_{n,\theta}(X_n)$ is given by (1.2) and $\alpha^*(n, k, \varepsilon) = \alpha$.

2.5. *Contamination tolerance of a test.* In some cases a test may be significant due to the presence of a small fraction of contamination in the data. To what extent might this be the case in a given application? The significance of the test would deliver a stronger message if we could discard the possibility that the results are due to contamination in the data. This motivates the following definition.

DEFINITION T3 (Contamination tolerance). Consider a family of tests $\varphi_{\theta_0,\varepsilon}$ for $H_0: \theta = \theta_0$ versus $H_1: \theta \neq \theta_0$, $0 \leq \varepsilon < 0.5$, such that (i) $\varphi_{\theta_0,\varepsilon}$ is $\varepsilon$-robust of level $\alpha$ and (ii) $\varepsilon_1 < \varepsilon_2$ implies $\varphi_{\theta_0,\varepsilon_1}(X_n) \geq \varphi_{\theta_0,\varepsilon_2}(X_n)$. Given a sample $X_n$ such that $\varphi_{\theta_0,0}(X_n) = 1$, the contamination tolerance for significance level $\alpha$ at $X_n$ [denoted by $\tau_\alpha = \tau_\alpha(X_n)$] is defined as

$$\tau_\alpha(X_n) = \sup\{\varepsilon : \varphi_{\theta_0,\varepsilon}(X_n) = 1\}.$$

In other words, the contamination tolerance for significance level $\alpha$ is the maximum level of contamination $\varepsilon$ such that the $\varepsilon$-robust test of level $\alpha$ still rejects the null hypothesis. Therefore, if we believe that the fraction of contamination in the data is smaller than $\tau_\alpha$, it is safe to reject the null hypothesis, even if we do not know the exact contamination size. Consequently, a large $\tau_\alpha$ (with small $\alpha$) can be taken as strong evidence against the null hypothesis.

Consider now the family of $\varepsilon$-robust sign tests given by (2.4). Then the value of $\tau_a(X_n)$ satisfies the equation

(2.5) $$\alpha^*\{n, r_n(X_n), \tau_a\} = \alpha,$$

where $r_n(X_n) = \min\{T_{n,\theta_0}(X_n), n - T_{n,\theta_0}(X_n)\}$. Notice that equation (2.5) has a solution if and only if $\alpha^*\{n, r_n(X_n), 0\} < \alpha$, that is, if and only if the null hypothesis is rejected under the assumption of a zero fraction of contamination (perfect data). If this condition is not satisfied, we would not reject $H_0$ even if the classical sign test is used.

**3. Length and power robustness.** Definitions I1 and I2 guarantee the correct coverage level of the interval. However, robust confidence intervals should not only have correct level but also remain informative under contamination. Definition I3 formalizes this robustness requirement in terms of the concept of *maximum asymptotic length* of the interval introduced below.

For the following discussion we must distinguish between the *design* contamination size $\varepsilon$ used to construct the confidence interval (so that it satisfies Definition I1) and the *real* contamination size denoted by $\delta$.

Given a sequence of intervals $I_n = [a_n(X_n), b_n(X_n))$, we consider the maximum asymptotic length under contamination of size $\delta$ at $F$,

(3.1) $$L\{I_n, F, \delta\} = \sup_{G \in \mathcal{F}_\delta(F)} \operatorname*{essup}\limsup_n (b_n(X_n) - a_n(X_n)),$$

where essup stand for essential supremum. The essup is applied for greater generality; however, in all cases we are aware of (including the interval based on the revised sign test), $\limsup_n (b_n(X_n) - a_n(X_n))$ is a constant (finite or infinite) and, therefore, essup is not necessary. Notice that if the interval length is location invariant, so is the above definition.

The intuitive notion of remaining "informative under contamination of size $\delta$" is captured by the following definition. Notice that our definition of length breakdown point is the confidence interval counterpart of Hampel's (1971) breakdown point of a point estimate.

DEFINITION I3 (Length robustness). We say that the sequence of intervals $I_n = [a_n(X_n), b_n(X_n))$, $n \geq n_0$, has $\delta$-robust length at $F$ if $L\{I_n, F, \delta\} < \infty$. The corresponding length breakdown point at $F$ is given by

$$\delta^*\{I_n, F\} = \sup\{\delta : L\{I_n, F, \delta\} < \infty\}.$$

The next theorem establishes the asymptotic length-robustness of the modified sign test interval.

THEOREM 2. *Suppose that $F$ is continuous and has a symmetric (around $\theta$) and unimodal density. Let $0 < \alpha < 1$ and $0 \leq \varepsilon < 1/2$ be fixed and consider the sequence of intervals $I_n = [x_{(k_n+1)}, x_{(n-k_n)})$, with $k_n$ given by (2.3). Then:*

1. *For $0 \leq \delta < (1-\varepsilon)/2$,*

$$L\{I_n, F, \delta\} = F^{-1}\left\{\frac{1+\varepsilon}{2(1-\delta)}\right\} - F^{-1}\left\{\frac{1-\varepsilon}{2(1-\delta)}\right\}.$$

2. $\delta^*\{(I_n), F\} = (1-\varepsilon)/2$.
3. *The sequence of intervals $I_n$ has $\varepsilon$-robust length if and only if $\varepsilon < 1/3$.*
4. *Let $I_n = [A_n(X_n), B_n(X_n))$ be a sequence of confidence intervals such that*

$$\inf_{G \in \mathcal{F}_\varepsilon(G_0)} P_G\{A_n(X_n) \leq G_0^{-1}(1/2) < B_n(X_n)\} = 1 - \alpha$$

*for any continuous distribution $G_0$. Suppose that $\lim_{n \to \infty} A_n(X_n) = A_0$, and $\lim_{n \to \infty} B_n(X_n) = B_0$ almost surely when the sample comes from $F$. Then $B_0 \geq F^{-1}((1+\varepsilon)/2)$ and $A_0 \leq F^{-1}((1-\varepsilon)/2)$.*

As one may have expected, the maximum asymptotic length of the sign-test-based intervals depends on the design and actual fractions of contamination, $\varepsilon$ and $\delta$. Finite maximum lengths are obtained provided $\delta < (1-\varepsilon)/2$. Therefore, length-breakdown point occurs when $\delta = (1-\varepsilon)/2$. Since the length-breakdown point $\delta^* = (1-\varepsilon)/2$ is a decreasing function of $\varepsilon$, there is a trade-off between the coverage-robustness and the length-robustness of the sign-test-based intervals. This naturally sets an upper bound of 1/3 on the possible choices of design-contamination fractions in practice. Part 4 shows that in the case of uncontaminated data (i.e., $\delta = 0$), our interval is efficient in that it has the smallest possible asymptotic length among all nonparametric $\varepsilon$-robust confidence intervals for the median, which upper and lower limits converge. Notice that convergence of the interval limits is a weak assumption satisfied by all known confidence intervals.

TABLE 2
*Coverage probability (CP) and expected length (EL) for robust confidence interval with approximate 95% coverage probability*

|       | $\varepsilon = 0$ |      | $\varepsilon = 0.05$ |      |      | $\varepsilon = 0.10$ |      |      |
|-------|-------|------|-------|------|------|-------|------|------|
| $n$   | CP    | ELU  | CP    | ELU  | ELC  | CP    | ELU  | ELC  |
| 20    | 0.959 | 1.22 | 0.954 | 1.22 | 1.3  | 0.938 | 1.24 | 2.52 |
| 40    | 0.962 | 0.84 | 0.952 | 0.83 | 0.89 | 0.960 | 0.97 | 1.13 |
| 60    | 0.948 | 0.64 | 0.961 | 0.72 | 0.76 | 0.955 | 0.81 | 0.92 |
| 80    | 0.943 | 0.54 | 0.949 | 0.60 | 0.64 | 0.955 | 0.73 | 0.84 |
| 100   | 0.943 | 0.48 | 0.941 | 0.53 | 0.56 | 0.957 | 0.69 | 0.78 |
| 200   | 0.944 | 0.34 | 0.947 | 0.42 | 0.44 | 0.949 | 0.55 | 0.61 |
| 500   | 0.946 | 0.22 | 0.947 | 0.31 | 0.32 | 0.952 | 0.44 | 0.50 |
| 1000  | 0.946 | 0.15 | 0.947 | 0.25 | 0.27 | 0.948 | 0.38 | 0.43 |
| 2000  | 0.948 | 0.11 | 0.949 | 0.22 | 0.23 | 0.950 | 0.34 | 0.39 |

3.1. *Numerical results.* We wrote a simple S-PLUS function, available on-line at http://hajek.stat.ubc.ca/~ruben/code1, which for a given sample $X_n$, significance level $\alpha$, and design contamination fraction $\varepsilon$, reports the integer $k_n$, the robust interval $[x_{(k_n+1)}, x_{(n-k_n)})$ and its exact minimum coverage probability, $1 - \alpha^*(n, k_n, \varepsilon)$. Using this function, we carried out a Monte Carlo simulation study to determine the increase in expected length for the robust nonparametric intervals $[x_{(k_n+1)}, x_{(n-k_n)})$ with $k_n$ given by (2.3).

We consider two approximate coverage probabilities, 95% (Table 2) and 90% (Table 3) and three contamination levels $\varepsilon = 0$, 0.05 and 0.10. The case $\varepsilon = 0$ corresponds to confidence intervals based on the classical sign-test. The tables display the exact infimum coverage probabilities (CP) and average lengths (EL). The average lengths of the robust confidence intervals are computed under two scenarios: uncontaminated (ELU) and contaminated samples (ELC). In the latter case, the fraction of contamination ($\delta$) equals the design contamination ($\varepsilon$). The contamination is placed at the least favorable location, which, as shown in the proof of Theorem 2, corresponds to $H = \delta_y$ in (1.1) with $y \to \pm\infty$. Naturally, the percent increase in average length is larger for larger samples sizes, when the effect of sampling variability is overcome by the effect of contamination bias. The average lengths are computed using 8000 replications.

In Table 4 we compare the asymptotic length of the nonparametric robust confidence intervals with the limiting length of the asymptotic parametric robust confidence intervals proposed by Huber (1968) and Fraiman, Yohai and Zamar (2001). The latter were proposed for a contamination neighborhood of the normal distribution and have limiting length equal to $2\Phi^{-1}[1/\{2(1-\varepsilon)\}]$, which is twice the maximum asymptotic bias of the median over the contamination neighborhood. We calculated the limiting lengths for both proposals under the normal model and under the least favorable contaminating distribution in $\mathcal{F}_\varepsilon(\Phi)$.

Notice that under Standard Normal, the nonparametric robust intervals have smaller expected length for all the considered values of $\varepsilon$. The expected lengths are practically equal for the least favorable contamination with a small advantage for the parametric interval.

3.2. *Power robustness of a test.* As in the case of confidence intervals, we must distinguish between the design contamination $\varepsilon$ used to construct the test and the actual contamination $\delta$. The following definition formalizes the concept of *robust power behavior* of a test under contamination of size $\delta$.

TABLE 3
*Coverage probability (CP) and expected length (EL) for robust confidence interval with approximate 90% coverage probability*

|      | $\varepsilon = 0$ |      | $\varepsilon = 0.05$ |      |      | $\varepsilon = 0.10$ |      |      |
|------|------|------|------|------|------|------|------|------|
| $n$  | CP   | ELU  | CP   | ELU  | ELC  | CP   | ELU  | ELC  |
| 20   | 0.885 | 0.89 | 0.876 | 0.90 | 0.96 | 0.938 | 1.20 | 2.40 |
| 40   | 0.919 | 0.70 | 0.904 | 0.70 | 0.74 | 0.922 | 0.83 | 0.95 |
| 60   | 0.908 | 0.55 | 0.883 | 0.55 | 0.58 | 0.923 | 0.72 | 0.82 |
| 80   | 0.907 | 0.47 | 0.918 | 0.54 | 0.57 | 0.891 | 0.60 | 0.68 |
| 100  | 0.911 | 0.43 | 0.912 | 0.48 | 0.51 | 0.904 | 0.58 | 0.66 |
| 200  | 0.896 | 0.29 | 0.908 | 0.36 | 0.39 | 0.912 | 0.49 | 0.56 |
| 500  | 0.902 | 0.19 | 0.895 | 0.27 | 0.28 | 0.904 | 0.40 | 0.45 |
| 1000 | 0.906 | 0.13 | 0.903 | 0.23 | 0.24 | 0.904 | 0.36 | 0.40 |
| 2000 | 0.897 | 0.09 | 0.899 | 0.20 | 0.21 | 0.900 | 0.32 | 0.36 |

TABLE 4
*Expected length of parametric (P) and nonparametric (NP) robust intervals*

|                 | $\varepsilon = 0.05$ |       | $\varepsilon = 0.10$ |       | $\varepsilon = 0.15$ |       | $\varepsilon = 0.20$ |       |
|-----------------|-------|-------|-------|-------|-------|-------|-------|-------|
| Distribution    | P     | NP    | P     | NP    | P     | NP    | P     | NP    |
| Standard Normal | 0.132 | 0.125 | 0.279 | 0.251 | 0.446 | 0.378 | 0.637 | 0.507 |
| Least Favorable | 0.132 | 0.132 | 0.279 | 0.282 | 0.446 | 0.458 | 0.637 | 0.674 |

DEFINITION T4 (Power robustness). Let $F$ be a fixed distribution satisfying (A1) with $\theta = \theta_0$ and let $F_\lambda(x) = F(x - \lambda)$. We say that a sequence of nonrandomized tests $\{\varphi_{n,\theta_0}\}$, $n \geq n_0$, has $\delta$-robust power (for $H_0$ versus $H_1$) at $F$ if there exists $K$ such

$$(3.2) \quad \inf_{G \in \mathcal{F}_\delta(F_\lambda)} \lim_{n \to \infty} P_G\{\varphi_{n,\theta_0}(X_n) = 1\} = 1 \quad \text{for all } |\lambda| > K.$$

This property ensures the consistency of the sequence of nonrandomized tests $\{\varphi_{n,\theta_0}\}$, uniformly over the neighborhood $\mathcal{F}_\varepsilon(F_\lambda)$, provided $\lambda = \theta - \theta_0$ is large enough. Definition T4 suggests the following measure of asymptotic power robustness of the sequence $\{\varphi_{n,\theta_0}(X_n)\}$ of tests, under contamination of size $\delta$.

DEFINITION T5 (Power distance). Let $F$ be a fixed distribution satisfying (A1) with $\theta = \theta_0$. The $\delta$-consistency distance of a sequence of tests $\varphi_{n,\theta_0}$, $n \geq n_0$, at $F$ denoted by $K\{\varphi_{n,\theta_0}, F, \delta\}$ is the infimum of the set of values $K$ for which (3.2) holds.

The concept of breakdown point of a test was first consider by Ylvisaker (1977) and Rieder (1982). The latter defined and computed the breakdown point of rank and M-tests. Our Definition T5 is closely related to the concept of power breakdown point of a test introduced by He, Simpson and Portnoy (1990). In fact, for a given $\theta \neq \theta_0$, the power breakdown point at $\theta$ is the value of $\delta$ such that $|\theta - \theta_0| = K\{(\varphi_{n,\theta_0}), F, \delta\}$.

Next we define a new concept of breakdown point for a test which does not depend on a particular value of $\theta$ and is directly associated with the definition of length breakdown point of a confidence interval given in Section 3.

DEFINITION T6 (Power breakdown). Let $F$ be a fixed distribution satisfying $\theta_0 = F^{-1}(1/2)$. The power breakdown point $\delta^*$ of the sequence of nonrandomized tests $\varphi_{n,\theta_0}$, $n \geq n_0$, at $F$ is the supremum of the set of values $\delta$ for which the sequence of tests is $\delta$-robust.

The power-robustness properties of the robustified sign test given by (2.4) are established in the next theorem. They are closely related to the length-robustness properties of the confidence intervals established in Theorem 2.

THEOREM 3. *Let $0 < \alpha < 1$ and $0 \leq \varepsilon < 1/2$ be fixed and consider the sequence of tests $\varphi_{\theta_0, n}$, $n \geq n_0$, for $H_0 : \theta = \theta_0$ versus $H_1 : \theta \neq \theta_0$ given by (2.4) and $k_n$ given by (2.3). Suppose that $F$ is continuous and has a symmetric (around $\theta$) and unimodal density. Then:*

1. *The $\delta$-consistency distance for the sequence of tests $\varphi_{n,\theta_0}$, $n \geq n_0$, is*

$$K\{(\varphi_{n,\theta_0}), F, \varepsilon\} = F^{-1}\left\{\frac{1+\varepsilon}{2(1-\delta)}\right\}.$$

2. *The power breakdown point of the sequence of tests $\varphi_{n,\theta_0}$, $n \geq n_0$, is $\delta^* = (1-\varepsilon)/2$.*
3. *The sequence of tests $\varphi_{n,\theta_0}$, $n \geq n_0$, has $\varepsilon$-robust power if and only if $\varepsilon < 1/3$.*

**4. Possible extensions and further research.** Robust nonparametric confidence intervals and tests for a location parameter could be defined using other rank statistics such as the signed Wilcoxon test statistics. In this case the parameter of interest would be defined as the center of symmetry of the target distribution, and, therefore, the target distribution (but not the observed distribution) would need to be symmetric. The main theoretical problem, which we were not able to solve, is the derivation of the least favorable distribution that gives the minimum coverage. We conjecture that this distribution is the one that puts all its mass at $+\infty$ or at $-\infty$.

We are currently studying possible extensions of our procedure to the case of two samples and to the case of simple linear regression. For the two-sample problem, we wish to construct robust nonparametric confidence intervals for the shift parameter, based on the two sample median test statistic. For the simple linear regression problem, we wish to construct robust nonparametric confidence intervals for the slope parameter, based on the Brown and Mood (1951) test statistic, which is a natural extension of the sign test statistic.

## APPENDIX

Lemma 1 is needed to prove Theorem 1. The proof of this lemma can be found as Lemma 4 of Yohai and Zamar (2004).

LEMMA 1. *Suppose that $X$ is $\mathrm{Bin}(n,p)$ and let*

$$h(p) = \sum_{i=k}^{n-k} \binom{n}{i} p^i (1-p)^{n-i}.$$

*Then* (i) $h(p) = h(1-p)$, (ii) $h(p)$ *is nondecreasing on $0 \leq p \leq 1/2$ for all $k = 0, 1, \ldots, [n/2]$.*

PROOF OF THEOREM 1. We have

$$P_G(x_{(k+1)} \leq \theta < x_{(n-k)}) = P_G\{k < T_{n,\theta}(X_n) < n - k\}$$
$$= P(k < Z_n < n - k),$$

where $Z_n$ is distributed as Binomial$\{n, 1-G(\theta)\}$. On the other hand, $G(\theta) = (1-\varepsilon)F(\theta) + \varepsilon H(\theta)$ and so

$$\frac{1-\varepsilon}{2} = (1-\varepsilon)F(\theta) \leq G(\theta) \leq (1-\varepsilon)F(\theta) + \varepsilon = \frac{1+\varepsilon}{2}.$$

Therefore, for all $G \in \mathcal{F}_\varepsilon(F), (1-\varepsilon)/2 \leq 1 - G(\theta) \leq (1+\varepsilon)/2$ with the lower and upper bounds attained when $H(\theta)$ concentrates all its mass to the left and right of $\theta$, respectively. The theorem now follows from Lemma 1. □

The following lemma is needed to prove Theorem 2. For a proof of Lemma 2 see Lemma 5 in Yohai and Zamar (2004).

LEMMA 2. *Let $X_n = (x_1, \ldots, x_n)$ be i.i.d. random variables with distribution $G$. Consider the sequence of intervals $I_n(X_n) = [x_{(k_n+1)}, x_{(n-k_n)})$ with lengths $l_n(X_n) = x_{(n-k_n)} - x_{(k_n+1)}$ and levels $\alpha^*(n, k_n, \varepsilon) \to \alpha$, $0 < \alpha < 1$. Then $\lim_{n \to \infty} l(X_n) = G^{-1}(\frac{1+\varepsilon}{2}) - G^{-1}(\frac{1-\varepsilon}{2}) = L(G, \varepsilon)$.*

PROOF OF THEOREM 2. Put $L^*(G, \varepsilon) = G^{-1}\{(1+\varepsilon)/2\} - G^{-1}\{(1-\varepsilon)/2\}$. By Lemma 2, to prove part 1 it is enough to show

$$(A.1) \qquad \sup_{G \in \mathcal{F}_\delta(F)} L^*(G, \varepsilon) = F^{-1}\left\{\frac{1+\varepsilon}{2(1-\delta)}\right\} - F^{-1}\left\{\frac{1-\varepsilon}{2(1-\delta)}\right\}.$$

We start by showing that

$$(A.2) \qquad \sup_{G \in \mathcal{F}_\delta(F)} L^*(G, \varepsilon) \leq F^{-1}\left\{\frac{1+\varepsilon}{2(1-\delta)}\right\} - F^{-1}\left\{\frac{1-\varepsilon}{2(1-\delta)}\right\}.$$

Let $G = (1-\delta)F + \delta H$. Then

$$a_1 = G^{-1}\left(\frac{1-\varepsilon}{2}\right), \qquad a_2 = G^{-1}\left(\frac{1+\varepsilon}{2}\right),$$

$$a_3 = F^{-1}\left\{\frac{1-\varepsilon}{2(1-\delta)}\right\} = 0, \qquad a_4 = F^{-1}\left\{\frac{1+\varepsilon}{2(1-\delta)}\right\}.$$

We will show first that

$$(A.3) \qquad F(a_2) - F(a_1) \leq F(a_4) - F(a_3).$$

This follows because by definition of quantiles,

$$\varepsilon = G(a_2) - G(a_1) = (1-\delta)F(a_2) + \delta H(a_2) - (1-\delta)F(a_1) - \delta H(a_1)$$
$$= (1-\delta)\{F(a_2) - F(a_1)\} + \delta\{H(a_2) - H(a_1)\}$$
$$\geq (1-\delta)\{F(a_2) - F(a_1)\},$$

and, therefore,

(A.4) $$F(a_2) - F(a_1) \leq \frac{\varepsilon}{1-\delta}.$$

On the other hand,

(A.5) $$F(a_4) - F(a_3) = \frac{1+\varepsilon}{2(1-\delta)} - \frac{1-\varepsilon}{2(1-\delta)} = \frac{\varepsilon}{1-\delta}.$$

Therefore, (A.3) follows from (A.4) and (A.5). To complete the proof of (A.2), we consider two cases:

*Case* 1 ($\delta \leq \varepsilon$). First notice that:

(i) $1/2 \geq (1-\varepsilon)/\{2(1-\delta)\}$ implies that

$$0 = F^{-1}\left(\frac{1}{2}\right) \geq F^{-1}\left\{\frac{1-\varepsilon}{2(1-\delta)}\right\} = a_3.$$

(ii) $(1-\varepsilon)/2 = F(a_1) \geq (1-\delta)F(a_1)$ implies that

$$a_1 \leq F^{-1}\left\{\frac{1-\varepsilon}{2(1-\delta)}\right\} = a_3.$$

By (A.3),

(A.6) $$F(a_4) - F(a_2) \geq F(a_3) - F(a_1).$$

Given the symmetry and unimodality of $F$, (A.2) follows from (A.6) if we can show that

(A.7) $$a_2 \geq -a_3.$$

To prove (A.7), we first notice the identity

(A.8) $$\frac{(\varepsilon-\delta)}{2(1-\delta)} = \frac{1}{2} - \frac{1-\varepsilon}{2(1-\delta)} = \frac{1+\varepsilon-2\delta}{2(1-\delta)} - \frac{1}{2}.$$

Symmetry of $F$ and (A.8) imply

(A.9) $$F^{-1}\left\{\frac{1+\varepsilon-2\delta}{2(1-\delta)}\right\} = -F^{-1}\left\{\frac{1-\varepsilon}{2(1-\delta)}\right\} = -a_3.$$

Moreover, $(1+\varepsilon)/2 = G(a_2) \leq (1-\delta)F(a_2) + \delta$ implies

(A.10) $$a_2 \geq F^{-1}\left\{\frac{1+\varepsilon-2\delta}{2(1-\delta)}\right\}.$$

Equation (A.7) follows now from (A.9) and (A.10).

*Case* 2 ($\delta > \varepsilon$). Since in this case $1/2 < (1-\varepsilon)/\{2(1-\delta)\}$, we have

(A.11) $$0 = F^{-1}\left(\frac{1}{2}\right) \leq F^{-1}\left\{\frac{1-\varepsilon}{2(1-\delta)}\right\} = a_3.$$

Moreover, $(1-\varepsilon)/2 = G(a_1) \leq (1-\delta)F(a_1) + \delta$ implies

(A.12) $$a_1 \geq F^{-1}\left\{\frac{1-\varepsilon-2\delta}{2(1-\delta)}\right\}.$$

We have the identity

(A.13) $$\frac{\varepsilon+\delta}{2(1-\delta)} = \frac{1}{2} - \frac{1-\varepsilon-2\delta}{2(1-\delta)} = \frac{1+\varepsilon}{2(1-\delta)} - \frac{1}{2}.$$

Equations (A.12) and (A.13) give

(A.14) $$a_1 \geq F^{-1}\left\{\frac{1-\varepsilon-2\delta}{2(1-\delta)}\right\} = -F^{-1}\left\{\frac{1+\varepsilon}{2(1-\delta)}\right\} = -a_4.$$

The inequality $(1-\varepsilon)/2 = G(a_1) \geq (1-\delta)F(a_1)$ implies

(A.15) $$a_1 \leq F^{-1}\left\{\frac{1-\varepsilon}{2(1-\delta)}\right\} = a_3.$$

Equations (A.14) and (A.15) give

(A.16) $$-a_4 \leq a_1 \leq a_3.$$

Then (A.2) follows now from (A.16) and the unimodality and the symmetry of $F$.

Let $\delta_m$ be the point mass distribution at $m$. Then
$$\lim_{m\to\infty} L^*\{(1-\varepsilon)F + \varepsilon\delta_m, \varepsilon\} = F^{-1}\left\{\frac{1+\varepsilon}{2(1-\delta)}\right\} - F^{-1}\left\{\frac{1-\varepsilon}{2(1-\delta)}\right\}.$$

This together with (A.2) implies (A.1). The proofs of parts 2 and 3 are straightforward. To prove part 2 just notice that the maximum interval length is finite provided that $(1+\varepsilon)/\{2(1-\delta)\} < 1$. Part 3 follows immediately from part 2.

Finally, to prove part 4, let $G_0$ be defined by

$$G_0(x) = \begin{cases} 0, & \text{if } x < F^{-1}(\varepsilon), \\ \dfrac{F(x)-\varepsilon}{1-\varepsilon}, & \text{if } x \geq F^{-1}(\varepsilon), \end{cases}$$

and $H$ be defined by

$$H(x) = \begin{cases} \dfrac{F(x)}{\varepsilon}, & \text{if } x < F^{-1}(\varepsilon), \\ 1, & \text{if } x \geq F^{-1}(\varepsilon). \end{cases}$$

Then observe that $F = (1-\varepsilon)G_0 + \varepsilon H$, and, therefore, $F \in \mathcal{F}_\varepsilon(G_0)$. Consequently, $G_0^{-1}(1/2) = F^{-1}((1+\varepsilon)/2) \in [A_0, B_0]$ and, therefore, $B_0 \geq F^{-1}((1+\varepsilon)/2)$.

Put

$$G_0(x) = \begin{cases} \dfrac{F(x)}{1-\varepsilon}, & \text{if } x < F^{-1}(1-\varepsilon), \\ 1, & \text{if } x \geq F^{-1}(1-\varepsilon), \end{cases}$$

$$H(x) = \begin{cases} 0, & \text{if } x < F^{-1}(1-\varepsilon), \\ \dfrac{F(x) - (1-\varepsilon)}{\varepsilon}, & \text{if } x \geq F^{-1}(1-\varepsilon). \end{cases}$$

We also have that $F = (1-\varepsilon)G_0 + \varepsilon H$, and, therefore, $F \in \mathcal{F}_\varepsilon(G_0)$. Then $G_0^{-1}(1/2) = F^{-1}((1-\varepsilon)/2) \in [a_0, b_0]$ and, therefore, $A_0 \leq F^{-1}((1-\varepsilon)/2)$. □

PROOF OF THEOREM 3. We can assume without loss of generality that $\theta_0 = 0$. We start by showing that given any $F$, we have

$$\text{(A.17)} \lim_{n \to \infty} P_G(\varphi_{n,0} = 1) = \begin{cases} 1, & \text{if } G^{-1}\{(1-\varepsilon)/2\} > 0, \\ 0, & \text{if } G^{-1}\{(1-\varepsilon)/2\} < 0 < G^{-1}\{(1+\varepsilon)/2\}, \\ 1, & \text{if } G^{-1}\{(1+\varepsilon)/2\} < 0. \end{cases}$$

We have

$$\text{(A.18)} \qquad P_G\{\varphi_{n,0}(X_n) = 1\} = P_G\{0 \notin [x_{(k_n)}, x_{(n-k_n)}]\}.$$

In Lemma 2 we have shown that $x_{(k_n)} \to G^{-1}\{(1-\varepsilon)/2\}$ and $x_{(n-k_n)} \to G^{-1}\{(1+\varepsilon)/2\}$. Therefore, (A.17) follows from (A.18). Then

$$\inf_{G \in \mathcal{F}_\delta(F_\lambda)} \lim_{n \to \infty} P_G\{\varphi_{n,0}(X_n) = 1\} = 1 \qquad \text{for all } |\lambda| > K$$

holds either if

$$\text{(A.19)} \qquad \sup_{G \in \mathcal{F}_\delta(F_\lambda)} G^{-1}\left(\frac{1+\varepsilon}{2}\right) = \lambda + \sup_{G \in \mathcal{F}_\delta(F)} G^{-1}\left(\frac{1+\varepsilon}{2}\right) < 0$$

or

$$\text{(A.20)} \qquad \inf_{G \in \mathcal{F}_\delta(F_\lambda)} G^{-1}\left(\frac{1-\varepsilon}{2}\right) = \lambda + \inf_{G \in \mathcal{F}_\delta(F)} G^{-1}\left(\frac{1-\varepsilon}{2}\right) > 0.$$

As in Theorem 2, we can show that

$$\sup_{G \in \mathcal{F}_\delta(F)} G^{-1}\left(\frac{1+\varepsilon}{2}\right) = F^{-1}\left\{\frac{1+\varepsilon}{2(1-\delta)}\right\}$$

and

$$\inf_{G \in \mathcal{F}_\delta(F)} G^{-1}\left(\frac{1-\varepsilon}{2}\right) = F^{-1}\left\{\frac{1-\varepsilon-2\delta}{2(1-\delta)}\right\} = -F^{-1}\left\{\frac{1+\varepsilon}{2(1-\delta)}\right\}.$$

In order for either (A.19) or (A.20) to hold, it is required that

$$|\lambda| > F^{-1}\left\{\frac{1+\varepsilon}{2(1-\delta)}\right\},$$

proving part 1 of the theorem. The proofs of parts 2 and 3 are straightforward. □

Departamento de Matemáticas  
Facultad de Ciencias Exactas y Naturales  
Universidad de Buenos Aires  
Ciudad Universitaria, Pabellón 1  
1428 Buenos Aires  
Argentina  
e-mail: vyohai@dm.uba.ar

Department of Statistics  
University of British Columbia  
6356 Agricultural Road  
Vancouver, British Columbia  
Canada V6T 1Z2  
e-mail: ruben@stat.ubc.ca